\documentclass{SCIYA2017enOL}

\def\a{\alpha} \def\b{\beta} \def\g{\gamma} \def\d{\delta} 
 \def\s{\sigma}   
 \def\ld{\lambda}

\def\G{\Gamma}

 \def\lg{\langle} \def\rg{\rangle}
\def\nd{\mathrel{\bigm|\kern-.7em/}}

\def\f{\noindent}

\def\PSL{\hbox{\rm PSL}}

\def\PGL{\hbox{\rm PGL}}

\def\Aut{\hbox{\rm Aut}}

\def\Cay{\hbox{\rm Cay}}

\def\Cos{\hbox{\rm Cos}}
\def\mod{\hbox{\rm mod }}

\def\demo{\f {\bf Proof.}\hskip10pt}

\def\mz{{\mathbb Z}}

\def\B{\mathcal{B}}
\def\F{\hbox{\rm F}}

\newtheorem{lem}{Lemma}[section]%
\newtheorem{defi}[lem]{Definition}%
\newtheorem{prop}[lem]{Proposition}%

\allowdisplaybreaks[4]
\online
\begin{document}

\ensubject{fdsfd}

\ArticleType{ARTICLES}
\Year{2017}
\Month{January}%
\Vol{60}
\No{1}
\BeginPage{1} %
\ReceiveDate{January 20, 2016}
\AcceptDate{May 5, 2017}

\title[]{Cubic vertex-transitive non-Cayley graphs \\ of order \textbf{12\textit{p}}}
{Cubic vertex-transitive non-Cayley graphs of order $12p$}

\author[1]{ZHANG Wei-Juan}{weijuanzhang369@gmail.com}
\author[2]{FENG Yan-Quan}{yqfeng@bjtu.edu.cn}
\author[2,$\ast$]{ZHOU Jin-Xin}{jxzhou@bjtu.edu.cn}

\AuthorMark{ZHANG}

\AuthorCitation{ZHANG W-J, FENG Y-Q, ZHOU J-X}

\address[1]{School of Mathematics (Zhuhai), Sun Yat-sen University, Guangdong {\rm 519082}, China }
\address[2]{Department of Mathematics, Beijing Jiaotong University, Beijing {\rm 100044}, China}

\abstract{A graph is said to be {\em vertex-transitive non-Cayley} if its full automorphism group acts transitively
on its vertices and contains no
subgroups acting regularly on its vertices. In this paper, a complete classification of cubic
vertex-transitive non-Cayley graphs of order $12p$, where $p$ is a prime, is given.
As a result, there are $11$ sporadic and one infinite family of such graphs, of which the sporadic ones occur when
$p=5$, $7$ or $17$, and the infinite family exists if and only if $p\equiv1\ (\mod 4)$,
and in this family there is a unique graph for a given order.}

\keywords{Cayley graphs, vertex-transitive graphs, automorphism groups}

\MSC{05C25, 20B25}

\maketitle


\section{Introduction}
Throughout this paper a {\em graph} means a finite, connected, simple and undirected graph.
Let $X$ be a graph with vertex set $V(X)$ and edge set $E(X)$.
Two vertices of $X$ are said to be {\em adjacent} if there is an edge between them.
An {\em arc} of $X$ is an ordered pair of adjacent vertices, and we use $A(X)$ to denote
the set of arcs of $X$. For $u,v\in V(X)$, $u\sim v$ means that $u$ is adjacent to $v$ and denote by $\{u,v\}$
the edge incident to $u$ and $v$ in $X$.
A graph $X$ is said to be {\em vertex-transitive} and
{\em arc-transitive} (or {\em symmetric}) if its full automorphism group, denoted by $\Aut(X)$,
acts transitively on $V(X)$ and $A(X)$, respectively.

Given a finite
group $G$ and a subset $S\subseteq G\setminus\{1\}$ such that $S=S^{-1}=\{s^{-1}\ |\ s\in S\}$ (where $1$
is the identity element of $G$),
the {\em Cayley graph} $\Cay(G,S)$ on $G$ with respect to $S$ is
defined to have vertex set $G$ and edge set $\{\{g,sg\}\mid g\in
G,s\in S\}$. A Cayley graph $\Cay(G,S)$ is connected if and only if
$S$ generates $G$. Given $g\in G$, define the permutation $R(g)$
on $G$ by $x\mapsto xg, x\in G$. Then $R(G)=\{R(g)\ |\ g\in G\}$,
called the {\em right regular representation} of $G$, is a
permutation group isomorphic to $G$. It is well-known that
$R(G)\leq\Aut(\Cay(G,S))$. So, $\Cay(G,S)$ is vertex-transitive.

A vertex-transitive graph is Cayley if
and only if its automorphism group contains a subgroup acting
regularly on its vertex set (see, for example, \cite[Lemma~4]{S0}).
This fact implies that Cayley graphs are just those vertex-transitive graphs
whose full automorphism groups have a regular subgroup.
Not all vertex-transitive graphs are Cayley and
the smallest one is the well-known Petersen graph. Such a graph will
be called a {\em vertex-transitive non-Cayley graph}, or a {\em
VNC-graph} for short.

In the literature, VNC-graphs have received considerable attention, and
much of the work have been focused on construction of VNC-graphs of certain type, determination of all VNC-graphs up to a certain order, classification of VNC-graphs with both specific valency and order, and so on.
For example, in 1983, Maru\v si\v c~\cite{M1} started a large project which aims at determining the set {\em NC} of non-Cayley numbers, that is, those numbers for which there exists a VNC-graph of order $n$. This was an active topic of research for a long time and
a lot of VNC-graphs were constructed in
\cite{HIP,LS,M2,MPr,MS,MS1,MSZ,MPr1,MPR,SE}. Motivated by the above mentioned research, Feng~\cite{F} considered the question to determine the smallest valency
$\vartheta(n)$ for
VNC-graphs of a given order $n$, and he answered this question for the graphs of odd prime power order.

In \cite[Table~1]{MPr}, the total number of
vertex-transitive graphs of order $n$ and the number of VNC-graphs
of order $n$ were listed for each $n\leq 26$, and very recently, Poto\v cnik et al.~\cite{PSV} constructed a census of all connected cubic vertex-transitive graphs on at most 1280 vertices. It seems that
Cayley graphs are `common' among the vertex-transitive graphs. In fact, McKay and Praeger conjectured that almost all vertex-transitive graphs are Cayley graphs.
This is true particularly for vertex-transitive graphs
with specific valency (see \cite{McKay-Royal,PSV}). These facts motivate us to study the problem of classifying VNC-graphs with specific valency.
From~\cite{CO,Marusic} we can obtain all VNC-graphs of order $2p$ for
each prime $p$. In \cite{Zhou-JGT-11} all tetravalent
VNC-graphs of order $4p$ were classified. In
\cite{Zhou-SSMS,Zhou-AdvM,ZF2010}, the cubic VNC-graphs of order
a product of three primes were classified. For the classification of
cubic VNC-graphs of order a  product of four primes,  Zhou and Feng \cite{ZF2012}
classified cubic VNC-graphs of order $8p$, and Kutnar et al.~\cite{KMZ} classified cubic
VNC graphs of order $4p^{2}$. In this paper
we classify all cubic VNC-graphs of order $12p$. As a result, there are 11 sporadic and one infinite
family of cubic VNC graphs of order $12p$, of which the
sporadic ones occur when $p=5$, $7$ or $17$, the infinite family exists if and only if
$p\equiv1\ (\mod 4)$, and in this family, there is one and only one graph for a given order.

\section{Preliminaries}
In this section, we introduce some notations and definitions as well
as some preliminary results which will be used later in the paper.
For a positive integer $n$, denote by $\mz_n$ the cyclic group of
order $n$ as well as the ring of integers modulo $n$, by $\mz_n^*$
the multiplicative group of $\mz_n$ consisting of numbers coprime to
$n$, by $D_{2n}$ the dihedral group of order $2n$, by $S_n$ the symmetric group of degree $n$,
and by $C_n$ and
$K_n$ the cycle and the complete graph of order $n$, respectively.
We call $C_n$ an {\em $n$-cycle}.

For two groups $M$ and $N$, $N\rtimes M$ denotes a semidirect
product of $N$ by $M$. For a subgroup $H$ of a group $G$, denote by
$C_G(H)$ the centralizer of $H$ in $G$ and by $N_G(H)$ the
normalizer of $H$ in $G$. Then $C_G(H)$ is normal in $N_G(H)$.
\begin{prop} {\rm\cite[Chapter I, Theorem~4.5]{H}}\ \  \label{NC}
The quotient group $N_G(H)/C_G(H)$ is isomorphic to a subgroup of
the automorphism group $\Aut(H)$ of $H$.
\end{prop}

Let $G$ be a permutation group on a set $\Omega$ and $\a\in \Omega$.
Denote by $G_\a$ the stabilizer of $\a$ in $G$, that is, the
subgroup of $G$ fixing the point $\a$. We say that $G$ is {\em
semiregular} on $\Omega$ if $G_\a=1$ for every $\a\in \Omega$ and
{\em regular} if $G$ is transitive and semiregular. For any $g\in
G$, $g$ is said to be {\em semiregular} if $\lg g\rg$ is
semiregular. The following proposition gives a characterization for
Cayley graphs in terms of their automorphism groups.

\begin{prop}{\rm\cite[Lemma~4]{S0}}\label{cayley graph}\ \ A graph $X$ is isomorphic
to a Cayley graph on a group $G$ if and only if its automorphism
group has a subgroup isomorphic to $G$, acting regularly on the
vertex set of $X$.
\end{prop}

An {\em $s$-arc} in a graph $X$ is an ordered $(s+1)$-tuple
$(v_0,v_1, \cdots ,v_{s-1},v_s)$ of vertices of $X$ such that
$v_{i-1}$ is adjacent to $v_i$ for $1\leq i\leq s$ and $v_{i-1}
\neq v_{i+1}$ for $1 \leq i \leq s-1$. A graph $X$ is said to be
{\em $s$-arc-transitive} if $\Aut(X)$ is transitive on the set of
$s$-arcs in $X$. In particular, $0$-arc-transitive means {\em
vertex-transitive}, and $1$-arc-transitive means {\em
arc-transitive} or {\em symmetric}. A subgroup of $\Aut(X)$ is
{\em $s$-regular} if the subgroup acts regularly on the set of
$s$-arcs in $X$, and $X$ is said to be {\em $s$-regular} if
$\Aut(X)$ is $s$-regular. Tutte \cite{T} proved that there exist
no cubic $s$-regular graphs for $s\geq 6$. The next proposition
characterizes the vertex stabilizers of cubic $s$-regular graphs.

\begin{prop} {\rm \label{prop3}}{\rm\cite[Propositions~2--5]{DM}}
For an $s$-regular cubic graph $X$, the stabilizer of $v\in V(X)$ in $\Aut(X)$ is isomorphic to $\mz_3$,
$S_3$, $S_3\times\mz_2$, $S_4$, or $S_4\times \mz_2$ for $s=1,2,3,4$
or $5$, respectively.
\end{prop}

For a regular graph $X$, use $d(X)$ to represent the valency of $X$,
and for any subset $B$ of $V(X)$, the subgraph of $X$ induced by $B$
will be denoted by $X[B]$. Let $X$ be a connected vertex-transitive
graph, and let $G\leq \Aut(X)$ be vertex-transitive on $X$. For a
$G$-invariant partition $\B$ of $V(X)$, the {\em quotient graph}
$X_\B$ is defined as the graph with vertex set $\B$ such that, for
any two vertices $B,C\in \B$, $B$ is adjacent to $C$ if and only if
there exist $u\in B$ and $v\in C$ which are adjacent in $X$. Let $N$
be a normal subgroup of $G$. Then the set $\B$ of orbits of $N$ in
$V(X)$ is a $G$-invariant partition of $V(X)$. In this case, the
symbol $X_\B$ will be replaced by $X_N$.

Let $X$ be a cubic graph and $G$ be an $s$-regular
subgroup of $\Aut(X)$ for some integer $s$.
Let $N$ be a normal subgroup of $G$.
In view of \cite[Theorem~9]{L1}, we have the following proposition.

\begin{prop} \label{Lorimer-th}
If $N$ has more than 
two orbits in $V(X)$, then $N$ is the kernel of $G$ on the set of orbits of $N$ that acts
semiregularly on $V(X)$, and $X_N$ is
a cubic symmetric graph with $G/N$ as an $s$-regular group of
automorphisms.
\end{prop}

A graph is said to be a {\em bi-Cayley graph} over a group $H$ if it admits $H$ as a semiregular automorphism
group with two orbits of equal size.
Note that every bi-Cayley graph admits the following concrete realization. Let $R,L$ and $S$ be subsets of a
group $H$ such that $R=R^{-1}$, $L=L^{-1}$ and $R\cup L$ does not contain the identity element of $H$. Define
the graph BiCay$(H, R, L, S)$ to have vertex set the union of the {\em right part} $H_0=\{h_0\ |\ h\in H\}$ and
the {\em left part} $H_1=\{h_1\ |\ h\in H\}$, and edge set the union of the {\em right edges}
$\{\{h_0,g_0\}\ |\ gh^{-1}\in R\}$, the {\em left edges} $\{\{h_1,g_1\}\ |\ gh^{-1}\in L\}$ and
the {\em spokes} $\{\{h_0,g_1\}\ |\ gh^{-1}\in S\}$.
For the case when $|S|=1$, the bi-Cayley graph BiCay$(H, R, L, S)$ is also called {\em one-matching bi-Cayley graph}
(see~\cite{kovacs}). Also, if $|R|=|L|=s$, then ${\rm BiCay}(H,R,L,S)$,
is said to be an {\em $s$-type bi-Cayley graph}.

To end this section, we introduce some properties of bi-Cayley graphs. In what follows, we always assume that $\G=$BiCay$(H,R,L,S)$ is a
connected bi-Cayley graph over a group $H$.

\begin{prop}{\rm\cite[Lemma~3.1]{ZF2014}}\label{basicprop}
The following hold.
\begin{enumerate}
  \item [{\rm (1)}]\  $H$ is generated by $R\cup L\cup S$.
  \item  [{\rm (2)}]\ $S$ can be chosen to contain the identity element of $H$.
  \item [{\rm (3)}]\ For any automorphism $\a$ of $H$, {\rm BiCay}$(H,R,L,S)\cong${\rm BiCay}$(H,R^\a,L^\a,S^\a)$.
\end{enumerate}
\end{prop}

Let $R(H)$ denote the right regular representation of $H$. Then
$R(H)$ can be regarded as a group of
automorphisms of BiCay$(H, R, L, S)$ acting on its vertices by the rule
$$h_i^{R(g)}=(hg)_i, \forall i\in\mz_2, h,g\in H.$$

For an automorphism $\a$ of $H$, define two permutations on $V(\G)=H_0\cup H_1$ as following:
\begin{equation}\label{eq-auto}
\begin{array}{ll}
\d_\a:& h_0\mapsto (h^\a)_1, h_1\mapsto (h^\a)_0, \forall h\in H,\\
\s_\a:& h_0\mapsto (h^\a)_0, h_1\mapsto (h^\a)_1, \forall h\in H.
\end{array}
\end{equation}
Set \begin{equation}\label{eq-autoset}
\begin{array}{lll}
{\rm I}&=& \{\d_\a\ |\ \a\in\Aut(H)\ s.t.\ R^\a=L, L^\a=R, S^\a=S^{-1}\},\\
{\rm F} &=&\lg \s_\a\ |\ \a\in\Aut(H)\ s.t.\ R^\a=R, L^\a=L, S^\a=S\rg.
\end{array}
\end{equation}

\begin{prop}{\rm \cite[Lemma~3.2 \& Theorem~3.4]{ZF2014}}\label{one-matching}
Let $\G=${\rm BiCay}$(H,R,L,\{1\})$ be a connected one-matching bi-Cayley graph over the group $H$.
Then $N_{\Aut(\G)}(R(H))=R(H)\rtimes {\rm F}$ if ${\rm I}=\emptyset$ and
$N_{\Aut(\G)}(R(H))=R(H)\rtimes \lg{\rm F},\d_\a\rg$ if ${\rm I}\not=\emptyset$ and $\d_\a\in{\rm I}$. Furthermore, if {\rm I} is non-empty, for any
$\d_\a\in{\rm I}$, we have the following:
\begin{enumerate}
  \item [{\rm(1)}]\ $\lg R(H),\d_\a\rg=R(H)\rtimes\lg \d_\a\rg$ acts transitively
on $V(\G)$;
  \item [{\rm (2)}]\ if $\a$ has order $2$, then $\G$ is isomorphic to the Cayley graph
  $\Cay(\bar{H}, R\cup \a S)$, where $\bar{H}= H\rtimes\lg \a\rg$.
\end{enumerate}
\end{prop}

\section{Cubic symmetric graphs of order $12p$}

In this section we shall classify cubic symmetric graphs of order $12p$ for each prime $p$.
Throughout this paper, the notations {\rm FnA}, {\rm FnB}, etc. will refer to the corresponding
graphs of order $n$ in the Foster census of all cubic symmetric graphs \cite{Bo,CD}.

\begin{theorem} {\rm \label{thm1}}
Let $X$ be a connected cubic symmetric graph of order $12p$ with $p$ a prime. Then $X$ is $2$- or $4$-regular, and moreover,
\begin{enumerate}
  \item [{\rm (1)}]\ $X$ is $2$-regular if and only if it is isomorphic to  $\F024$, $\F060$ or $\F084$;
  \item [{\rm (2)}]\ $X$ is $4$-regular if and only if it is isomorphic to $\F204$.
\end{enumerate}
\end{theorem}

\demo We first claim $p\leq 17$. Suppose to the contrary that $p>17$. Let $A=\Aut(X)$.
By Proposition~\ref{prop3}, the stabilizer $A_v$ of $v\in V(X)$ in $A$ has order dividing $48$.
It follows that $|A|=|V(X)||A_v|\ |\ 2^6\cdot3^2\cdot p$. If $A$ is non-solvable, then
$A$ has a non-abelian simple composite factor $M/N$, and since $|M/N|\ |\ 2^6\cdot3^2\cdot p$,
by \cite[pp.12--14]{G} $M/N$ is one of the following groups:
\begin{equation*}
A_{5},A_{6},\PSL(2,7),\PSL(2,8),\PSL(2,17),\PSL(3,3),{\rm PSU}(3,3),{\rm PSU}(4,2).
\end{equation*}
By considering the orders of these groups, we obtain that $p\leq 17$, a contradiction. Thus, $A$ is solvable.
Suppose that $A$ has a normal $r$-subgroup, say $N$, with $r=3$ or $p$.
Since $X$ has order $12p$ with $p>17$, it is easy to see that
$N$ has more than two orbits on $V(X)$. By Proposition~\ref{Lorimer-th},
$N$ is semiregular, and the quotient graph $X_N$ of $X$ relative to
$N$ is a cubic symmetric graph of order $12p/|N|$.
Since $N$ is semiregular, one has $|N|=3$ or $p$, implying that the number of vertices of $X_N$ is
$4p$ or $12$. However, by \cite[Theorem 6.2]{FK}, the cubic symmetric graph with order $12$
or $4p$ with $p>17$ does not exist, a contradiction. Thus, $A$ has no normal $3$- or $p$-subgroups.
Take a minimal normal subgroup, say $M$, of $A$. The solvability of $A$ implies that $M$
is an elementary abelian $r$-group, where $r\in\{2,3,p\}$. Then $M$ must be a $2$-group.
Clearly, $M$ has more than two orbits on $V(X)$. Again, by Proposition~\ref{Lorimer-th},
$M$ is semiregular, and the quotient graph $X_M$ of $X$ relative to
$M$ is a symmetric cubic graph. It follows that $M\cong\mz_2$, and $X_M$
has order $6p$. Let $T/M$ be a minimal normal subgroup of $A/M$.
Clearly, $A/M$ is solvable, so $T/M$ is an elementary abelian $r$-group, where $r\in\{2,3,p\}$.
If $T/M$ is a $r$-group with $r=3$ or $p$, then $T=R\times M$, where $R$ is the Sylow $r$-subgroup of
$T$. Clearly, $R$ is characteristic in $T$. The normality of $T$ in $A$ implies that
$R\unlhd A$. This is contrary to the fact that $A$ has no normal $3$- or $p$-subgroups.
If $T/M$ is a $2$-group, then $|T|=2^\ell\geq 4$. Clearly, $T$ has more than two orbits on $V(X)$.
By Proposition~\ref{Lorimer-th}, $T$ is semiregular and the quotient graph $X_T$ of $X$ relative to
$T$ is a cubic symmetric graph of odd order $12p/|T|=3p$, a contradiction.

By now, we have shown that our claim is true, namely, $p\leq 17$. Then
$X$ has order at most $204$, and by \cite{CD}, $X$ is isomorphic either to
the $2$-regular graphs $\F024, \F060$ or $\F084$, or to the $4$-regular graph $\F204$.\hfill\qed

\section{Graphs associated with lexicographic products}

Let $n$ be a positive integer. The {\em lexicographic product} $C_{n}[2K_1]$ is
defined as the graph with vertex set $\{x_i,y_i\ |\ i\in\mz_n\}$ and
edge set
$\{\{x_i,x_{i+1}\},\{y_i,y_{i+1}\},\{x_i,y_{i+1}\},\{y_i,x_{i+1}\}\
|\ i\in\mz_n\}$. In this section, we introduce a class of cubic vertex-transitive graphs which can be
constructed from the lexicographic product $C_n[2K_1]$. Note that these graphs
belong to a large family of graphs constructed in \cite[Section~3]{Dobson}.

\begin{defi} For integer $n\geqslant2$, let $X(n,2)$ be the graph of order $4n$
and valency $3$ with vertex set $V_{0}\cup V_{1}\cup\cdots V_{2n-2}\cup V_{2n-1}$, where
$V_{k}=\{x_{k}^{0},x_{k}^{1}\}$ $(k\in \mz_{2n})$, and adjacencies $x_{2i}^{r}\thicksim x_{2i+1}^{r}$
$(i\in \mz_{n},r\in \mz_{2})$ and $x_{2i+1}^{r}\thicksim x_{2i+2}^{s}$
$(i\in \mz_{n};r,s\in \mz_{2})$.\medskip
\end{defi}

Note that $X(n,2)$ is obtained from $C_n[2K_1]$ by expanding each vertex into an edge, in
a natural way, so that each of the two blown-up endvertices inherits half of the neighbors of the
original vertex.

Define three permutations on $V(X(n,2))$ as follows:
$$
\begin{array}{llll}
\a: & x_{k}^{0}\mapsto x_{k+2}^{0}, & x_{k}^{1}\mapsto x_{k+2}^{1},  & k\in\mz_{2n},\\
\b: &  x_{k}^{0}\mapsto x_{k}^{1},  & x_{k}^{1}\mapsto x_{k}^{0}, & k\in\mz_{2n},\\
\g: &  x_{k}^{0}\mapsto x_{2n-1-k}^{0},& x_{k}^{1}\mapsto x_{2n-1-k}^{1},& k\in\mz_{2n}.
\end{array}
$$
It is easy to check that $\alpha$, $\beta$, and $\gamma$
are automorphisms of $X(n,2)$. Furthermore, $\langle\alpha,\beta,\gamma\rangle\cong D_{4n}$ acts regularly
on $V(X(n,2))$, which implies that $X(n,2)$ is a Cayley graph.


\begin{theorem}\label{lem1}
Let $p>7$ be a prime, and $X$ a connected cubic vertex-transitive graph of order $12p$.
If $\Aut(X)$ is solvable, then either $\Aut(X)$ has a normal Sylow $p$-subgroup, or $X\cong X(3p,2)$.
\end{theorem}
\demo Let $A$=$\Aut(X)$. Since $X$ is cubic, the stabilizer $A_{v}$ of $v\in V(X)$
in $A$ has order dividing $2^m\cdot 3$ for some integer $m$, implying that
$|A|=|V(X)||A_{v}|\ |\ 2^{m+2}\cdot3^{2}\cdot p$. Let $P$ be a Sylow $p$-subgroup of $A$. Then $P\cong\mz_p$.
Denote by $O_r(A)$ the maximal normal $r$-subgroup of $A$ for some prime factor $r$
of $|A|$. If $O_{p}(A)>1$, then $P=O_p(A)$ is normal in $A$.
In what follows, assume that $O_p(A)=1$. We consider two cases: $O_2(A)=1$ and $O_2(A)>1$.\medskip

\f{\bf Case 1:}\ $O_{2}(A)=1$

In this case, we must have $O_{3}(A)>1$ because $A$ is solvable. Since $|A|\ |\ 2^{m+2}\cdot3^{2}\cdot p$,
one has $O_{3}(A)\cong \mz_{3}$, $\mz_{9}$ or $\mz_{3}^{2}$. Let $C:=C_{A}(O_{3}(A))$. Clearly, $O_{3}(A)\leq C$.
If $O_{3}(A)=C$, then
$A/O_3(A)=A/C\leq\Aut(O_{3}(A))$. Consequently, $p\ |\ |\Aut(O_{3}(A))|$. However,
$\Aut(O_{3}(A))\cong \mz_2, \mz_6$ or ${\rm GL}(2,3)$, whose orders can not be divided by $p$, a contradiction. Thus, $O_{3}(A)<C$.
Let $H/O_{3}(A)$ be a minimal normal subgroup of $A/O_{3}(A)$ such that $H/O_{3}(A)\leq C/O_{3}(A)$.
Then $H/O_{3}(A)$ is a $q$-group, where $q\in\{2,p\}$. Take a Sylow $q$-subgroup of $H$, say $Q$.
Since $Q\leq C$, one has $H=O_{3}(A)\times Q$, and hence $Q$ is characteristic in $H$.
The normality of $H$ in $A$ gives that
$Q\unlhd A$. This is contrary to the assumption that $O_2(A)=O_p(A)=1$.\medskip

\f {\bf Case 2:}\ $O_{2}(A)>1$

Let $X_{O_{2}(A)}$ be the quotient graph of $X$ relative to $O_{2}(A)$ with vertex set $\Omega$.
Remember that $\Omega$ is the set of the orbits of $O_2(A)$ in $V(X)$.
Let $K$ be the kernel of $A$ acting on $\Omega$.
Then $O_{2}(A)\leq K$, $A/K\leq \Aut(X_{O_{2}(A)})$
and $A/K$ acts transitively on $\Omega$.
In view of the fact that $O_{2}(A)>1$ and $X$ is of
order $12p$, we get that $X_{O_{2}(A)}$ has order $3p$ or $6p$. In addition,
since $X$ is cubic, the vertex-transitivity and connectivity of $X$ implies that the valency $d(X_{O_{2}(A)})$ of $X_{O_{2}(A)}$ is
$2$ or $3$.

If $d(X_{O_{2}(A)})=2$, then $X_{O_{2}(A)}$ is an $\ell$-cycle, where $\ell=3p$ or $6p$. Without loss
of generality, we may assume that
$\Omega=\{\Delta_{i}\ |\ i\in\mz_\ell\}$,
and $\Delta_{i}\sim\Delta_{i+1}$ for every $i\in \mz_{\ell}$.

Suppose first that $X[\Delta_{i}]$ is not an empty graph for an orbit $\Delta_{i}\in\Omega$.
Then, by the normality of $O_{2}(A)$, there is an edge inside each orbit in $\Omega$.
The fact that $X$ is cubic implies that the graphs induced by the orbits of $O_{2}(A)$ are either
all isomorphic to $2K_{2}$ or all isomorphic to $K_{2}$. Hence, for each $v\in \Delta_{i}$,
one neighbor of $v$ is in $\Delta_{i}$ and the other two neighbors are in $\Delta_{i-1}$ and $\Delta_{i+1}$,
respectively. Since $K$ fixes each orbit of $O_{2}(A)$ setwise, the stabilizer $K_{v}$ of $v\in V(X)$
fixes all neighbors of $v$. Applying the connectedness of $X$, we get $K_{v}$ fixes all vertices of $X$,
and hence $K_{v}=1$. It follows that $K=O_{2}(A)$ is semiregular, and so $|K|=12p/\ell=2$ or $4$.
Since $X_{O_{2}(A)}\cong C_{\ell}$,
one has $\Aut(X_{O_{2}(A)})\cong D_{2\ell}$. Since $A/K\leq \Aut(X_{O_{2}(A)})$ is transitive on
$\Omega$, one has $A/K\cong\mz_{\ell},D_{\ell}($for $\ell=6p)$ or $D_{2\ell}$.
Clearly, the Sylow $p$-subgroup $PK/K$ of $A/K$ is normal, so $PK\unlhd A$. Since $p>7$ and $|K|\leqslant4$,
Sylow Theorem yields $P\unlhd PK$, and hence $P$ is characteristic in $PK$.
The normality of $PK$ in $A$ gives that $P\unlhd A$, a contradiction.

Suppose now that $X[\Delta_{i}]$ is null for every $\Delta_{i}\in\Omega$. Since $X$ is cubic
and $X_{O_{2}(A)}$ is a cycle, the subgraph $X[\Delta_{i}\cup \Delta_{i+1}]$ of $X$ induced by any two
adjacent orbits $\Delta_{i}$ and $\Delta_{i+1}$ of $O_{2}(A)$ is a regular graph
with valency $1$ or $2$. Without loss of generality,
assume $d(X[\Delta_{0}\cup \Delta_{1}])=1$. Then $d(X[\Delta_{1}\cup \Delta_{2}])=2$. By vertex-transitivity of
$X$, we have $d(X[\Delta_{i}\cup \Delta_{i+1}])=1$ if $i$ is even, and $d(X[\Delta_{i}\cup \Delta_{i+1}])=2$
if $i$ is odd. Since $d(X[\Delta_{\ell-1}\cup \Delta_{0}])=2$, we must have $\ell$ is even. It follows that
$X_{O_{2}(A)}\cong C_{6p}$ and $|\Delta_{i}|=2$ for each $i\in \mz_{6p}$. Furthermore,
$X[\Delta_{i}\cup \Delta_{i+1}]\cong2K_{2}$ if $i$ is even whereas $X[\Delta_{i}\cup \Delta_{i+1}]\cong C_{4}$
if $i$ is odd. Let $\Delta_{i}=\{x_{i},y_{i}\}$, $i\in \mz_{6p}$. Then we may, without loss of generality,
assume that
$$E(X[\Delta_{i}\cup \Delta_{i+1}])=\left\{
\begin{array}{ll}
\{\{x_{i}, x_{i+1}\}, \{y_{i}, y_{i+1}\}\} & $~if~$i$~is~even$,\\
\{\{x_{i}, x_{i+1}\}, \{y_{i}, y_{i+1}\}, \{x_{i}, y_{i+1}\}, \{y_{i}, x_{i+1}\}\} &~$if~$i$~is~odd$. \\
\end{array}
\right.$$
Thus $X\cong X(3p,2)$.

If $d(X_{O_{2}(A)})=3$, then the length of every orbit of $O_{2}(A)$ is 2, and $X_{O_{2}(A)}$
has order $6p$.

Suppose that $X_{O_{2}(A)}$ is a non-Cayley graph. Since $p>7$,
by \cite[Theorem 5.3]{ZF2010}, $X_{O_{2}(A)}$ is isomorphic to the symmetric graph $\F102$
of order 102. Furthermore, the full automorphism group of $\F102$ is isomorphic to $\PSL(2,17)$.
Note that $A/K$ is a solvable vertex-transitive group of automorphisms of $X_{O_{2}(A)}$.
It follows that $A/K$ is a proper subgroup of $\PSL(2,17)$. Let $\tilde{M}$ be a maximal subgroup of
$\PSL(2,17)$ containing $A/K$. Then $102\ |\ |\tilde{M}|$ because $102\ |\ |A/K|$.
However, from \cite[Section 239]{Dickson}, we see that
$\tilde{M}$ is isomorphic to $D_{16}$, $D_{18}$, $\mz_{17}\rtimes \mz_{8}$ or $S_{4}$,
and none of their orders can be
divided by 102, a contradiction.

Suppose that $X_{O_{2}(A)}$ is a Cayley graph, say $X_{O_{2}(A)}\cong \Cay(G,S)$, where $G$ is a group of order $6p$.
Since $p>7$, by \cite[Theorem 3.2]{ZZF}, one has either $\Aut(X_{O_{2}(A)})\cong G\rtimes H$, where $H\leq S_3$, or
$\Aut(X_{O_{2}(A)})\cong G\mz_{3}$ for $3\ |\ (p-1)$.
For the former case, since $p>7$, the Sylow Theorem implies that the Sylow $p$-subgroup, say $R$, of $G$ is normal
and since $G\unlhd G\rtimes H$,
$R$ is also a normal Sylow $p$-subgroup of $G\rtimes H$. For the latter case, since $p>7$ and $3\ |\ (p-1)$,
again by Sylow Theorem, we get that the Sylow $p$-subgroup of $G\mz_{3}$ is normal.
As a result, $\Aut(X_{O_{2}(A)})$ always has a normal Sylow $p$-subgroup.
Since $A/K\leq \Aut(X_{O_{2}(A)})$, the Sylow $p$-subgroup $PK/K$ of $A/K$ is also normal,
and so $PK\unlhd A$. For every $v\in V(X)$, the neighbors of $v$ are in three different orbits of $O_{2}(A)$.
Since $K$ fixes each orbit of $O_{2}(A)$ setwise, $K_{v}$ fixes all neighbors of $v$. By the connectedness of
$X$, we have $K_{v}=1$ and $|O_{2}(A)|=2|K_{v}|=2$. Thus, $K=O_{2}(A)K_{v}=O_{2}(A)\cong \mz_{2}$ and
so $P\unlhd PK$. Since $P$ is a Sylow $p$-subgroup of $PK$, $P$ is characteristic in $PK$. The normality of
$PK$ in $A$ implies that $P\unlhd A$, contrary to the assumption of $O_p(A)=1$.\hfill\qed

\section{Construction of cubic VNC-graphs of order $12p$}

We shall first construct $9$ sporadic VNC-graphs of order $12p$ with $p=5$ or $7$.
To do so, we need the concept of coset graph (see \cite{S1}). The {\em coset graph}, denoted by $\Cos(G,H,D)$,
is constructed from a finite group $G$ relative to a subgroup $H$ of $G$ and a union $D$ of some double cosets
of $H$ in $G$ such that $D^{-1}=D$. The coset graph $\Cos(G,H,D)$ is defined to have vertex set $[G:H]$, the
set of right cosets of $H$ in $G$, and edge set $\{\{Hg, Hdg\}~|~g\in G, d\in D\}$. The valency of $\Cos(G,H,D)$
is $|D|/|H|=|H:H\cap H^{g}|$, and $\Cos(G,H,D)$ is connected if and only if $D$ generates the group $G$.
The action of $G$ on $V(\Cos(G,H,D))$ by right multiplication induces a vertex-transitive automorphism group.
$\Cos(G,H,D)$ is $G$-arc transitive if and only if $D$ is a single double coset.
In what follows, we shall construct 9 sporadic VNC-graphs of order $12p$ with $p=5$ or $7$ in term of coset graph.

\begin{defi}\label{def-1}
For $0\leq i\leq 8$, let ${\rm NC}_{12p}^i=\Cos(G,H,HaH\cup HbH)$, where $p=5$ or $7$, and $G,H,a,b$ are given in {\rm Table~\ref{k4t1}}.
\end{defi}
\begin{table}[ht]
\begin{center}
\begin{tabular}{|c|c|c|c|c|}
\hline
 Graph & $G$& $H$ & $a$ & $b$\\ \hline
$\rm NC_{12\cdot 5}^{0}$ &$S_{5}$&$\langle(1~3)(2~5)\rangle$&$(1~2~4~5~3)$&$(2~5)$ \\ \hline
$\rm NC_{12\cdot 5}^{1}$ &&&$(2~5)(3~4)(6~7)$&$(1~2)(3~4)(6~7)$ \\  \cline{4-5}\cline{1-1}
$\rm NC_{12\cdot 5}^{2}$ &$A_{5}\times \mz_{2}$&$\langle(2~3)(4~5)\rangle$&$(1~2~4~5~3)$&$(2~4)(3~5)(6~7)$ \\ \cline{4-5} \cline{1-1}
$\rm NC_{12\cdot 5}^{3}$ &&&$(1~3~5~4~2)(6~7)$&$(2~4)(3~5)$ \\ \hline
$\rm NC_{12\cdot 5}^{4}$ &&$\langle(1~3),(4~5)\rangle$&$(1~3~2)(4~5)$&$(1~5~3~4)(6~7)$ \\ \cline{3-5} \cline{1-1}
$\rm NC_{12\cdot 5}^{5}$ &$S_{5}\times \mz_{2}$&&$(2~5)$&$(1~4~3~5)(6~7)$ \\ \cline{4-5} \cline{1-1}
$\rm NC_{12\cdot 5}^{6}$ &&\raisebox{1.4ex}[0pt]{$\langle(1~3)(6~7),(4~5)(6~7)\rangle$}&$(2~4)$&$(1~4)(3~5)(6~7)$ \\ \hline
&&$\langle(1~7)(2~4)(3~8),$&& \\
\raisebox{1.4ex}[0pt]{$\rm NC_{12\cdot 7}^{7}$} &&$(1~4)(2~7)(5~6)\rangle$& \raisebox{1.4ex}[0pt]{$(1~8~4)(2~7~3)$}&
\raisebox{1.4ex}[0pt]{$(1~4~2~7)(3~6~8~5)$}\\ \cline{3-5}\cline{1-1}
&\raisebox{1.4ex}[0pt]{$\PGL(2,7)$}&$\langle(1~6)(2~4)(3~7)(5~8),$&& \\
\raisebox{1.4ex}[0pt]{$\rm NC_{12\cdot 7}^{8}$}&&$(1~4)(2~6)(3~8)(5~7)\rangle$&\raisebox{1.4ex}[0pt]{$(1~3~5~2~6~7~8~4)$}&
\raisebox{1.4ex}[0pt]{$(1~5~2~3)(4~8~6~7)$}
\\ \hline
\end{tabular}
\end{center}

\vskip -0.3cm
\caption{The VNC graphs of order 60\ or\ 84} \label{k4t1}
\end{table}

Using the computer software MAGMA~\cite{BCP}, we can easily obtain the following theorem.

\begin{theorem} {\rm \label{thm4}}
The graphs listed in Table~\ref{k4t1} are pairwise non-isomorphic non-symmetric cubic {\rm VNC}-graphs
with non-solvable automorphism groups.
\end{theorem}


Below, we shall construct an infinite family of VNC-graphs of order $12p$.
Let $p$ be an odd prime. It is well known that $\mz_{p}^*\cong\mz_{p-1}$.
So, if $4\ |\ (p-1)$ then $\mz_{p}^*$ has a unique subgroup of order $4$. Clearly, if $\ld$ is
an element of order $4$ in $\mz_{p}^*$, then $\{1,-1,\ld,-\ld\}$ is the unique subgroup of order $4$
in the cyclic group $\mz_{p}^*$.

\begin{defi}
Let $p$ be a prime congruent to $1$ modulo $4$ and let $\lambda$ be an element
of order $4$ in $\mz_p^*$. Let $H=\lg a,b\ |\ a^2=b^2=(ab)^3=1\rg\times\lg c\rg\cong S_3\times\mz_p$.
Set ${\rm NC}_{12p}^{9}:={\rm BiCay}(H, R, L, \{1\})$, where $R=\{ac, ac^{-1}\}$ and $L=\{bc^\ld, bc^{-\ld}\}$. See Figure 1 for the smallest one in this family of graphs.
\end{defi}

\begin{center}
\includegraphics[scale=0.6]{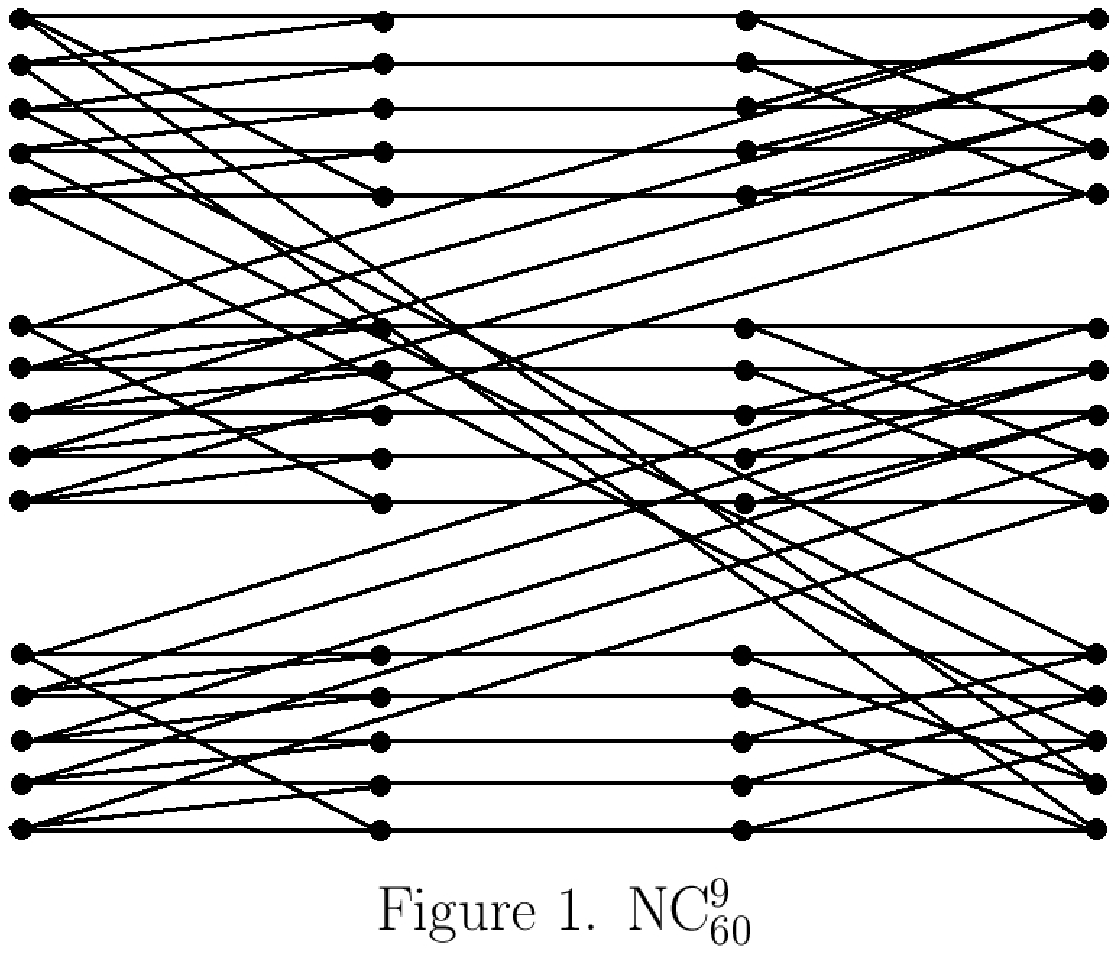}
\end{center}

The following theorem shows that
${\rm NC}_{12p}^{9}$ is a $\rm VNC$ graph of order $12p$.

\begin{theorem} {\rm \label{thm2}}
Let $p$ be a prime congruent to $1$ modular $4$. Then ${\rm NC}_{12p}^{9}$ is a connected cubic non-symmetric $\rm VNC$ graph of order $12p$,
and $\Aut({\rm NC}_{12p}^{9})\cong(S_{3}\times \mz_{p})\rtimes \mz_{4}$.
\end{theorem}
\demo Let $X={\rm NC}_{12p}^{9}$ and $A=$Aut$(X)$.
Clearly, $X$ is connected, and has valency three and order $12p$.
It is easy to see that $H$ has an automorphism, say $\a$ such that $a^\a=b, b^\a=a, c^\a=c^\ld$, and moreover, $\a$ interchanges $R$ and $L$. By Proposition~\ref{one-matching}, $\d_\a$ (see~(\ref{eq-auto})) is an automorphism of $X$ of order $4$ such that $\d_\a$ interchanges the two copies of $H$, and furthermore, $\d_\a\in N_A(R(H))$. Set $G=R(H)\rtimes\lg\d_\a\rg$. Then $G$ is vertex-transitive on $X$. Clearly, $G\leq A$. In what follows, we shall show that $G=A$.

If $p\leq 17$, then by MAGMA~\cite{BCP}, $X$ is non-symmetric, and by Theorem~\ref{thm1},
$X$ is also non-symmetric for $p>17$. So, the
stabilizer $A_{v}$ of $v\in V(X)$ in $A$ is a $2$-group, and hence
$|A|=|V(X)||A_{v}|\ |\ 2^{m+2}\cdot3\cdot p$ for a non-negative integer $m$.
As $p\geqslant5$, the group $P=\langle R(c)\rangle$ is a
Sylow $p$-subgroup of $A$. Again by MAGMA \cite{BCP}, if $p\leq 17$ then $A$ is solvable and contains a normal Sylow $p$-subgroup,
namely, $P=\lg R(c)\rg$.
For $p>17$, we have $X\ncong X(3p,2)$ because it is easy to check that $X$ has no $4$-cycles while
$X(3p,2)$ has girth $4$. It follows from Lemma \ref{lem1} that the
Sylow $p$-subgroup $P$ of $A$ is also normal.
Now consider the quotient graph $X_P$ of $X$ relative to $P$,
and let $K$ be the kernel of $A$ acting on $V(X_{P})$. From the construction of $X$, we see
$X_{P}\cong C_{12}$ and the subgraph of $X$ induced by any two adjacent orbits of $P$ is either isomorphic to $pK_{2}$
or $C_{2p}$. This implies that $K$ acts faithfully on each orbit of $P$, and hence $K\leq\Aut(C_{2p})\cong D_{4p}$.
Since $K$ fixes each orbit of $P$, one has $K\leq D_{2p}$. Clearly, $A/K$ is not edge-transitive on
$X_{P}$. It follows that $A/K\cong D_{12}$, and hence $|A|\leqslant24p$. Clearly,
$|G|=24p$, so $A=G$. 

Now we are ready to finish the proof.  Since $S_3\rtimes\mz_4=S_3\times \mz_4$, we have  $A=R(H)\rtimes\lg\d_\a\rg\cong (S_{3}\times \mz_{p})\rtimes \mz_{4}\cong S_3\times
(\mz_p\rtimes\mz_4)$, which has Sylow $2$-subgroups isomorphic to $\mz_2\times\mz_4$. Since Sylow $2$-subgroups are conjugate and $\d_\a$ interchanges the two orbits of $R(H)$ with $\d_\a^2$ fixing some vertices, every element of order $2$ in $A$ fixes the two orbits of $R(H)$ and every element of order $4$ interchanges the two orbits of $R(H)$ with its square fixing some vertices.
Suppose that $X$ is a Cayley graph. By Proposition \ref{cayley graph},
$A$ contains a regular subgroup $L$. Then $L$ contains a $2$-element $t$ of order $4$ interchanging the two orbits of $R(H)$. But $t^2$ fixes some vertices, contrary to the regularity of $L$.

Now we are ready to finish the proof. Suppose that $X$ is a Cayley graph.
By Proposition \ref{cayley graph},
$A$ contains a regular subgroup, say $U$. By $|A:U|=2$, we have $U$ is a maximal subgroup of $A$. Let $\g$ be an involution in $H$ and $Q=\lg \g, \d_\a\rg$. Then $Q$ is a Sylow 2-subgroup of $A$. Clearly, $Q\nleq U$ and thus $A=UQ$. It follows that $|Q\cap U|=4$ and $Q/(Q\cap U)\cong \mz_2$, and hence $\d_\a^2\in Q\cap U$. However, $a^{\d_\a^2}=a$, which implies that some vertices of $X$ are fixed by $\d_\a^2$, a contradiction. \hfill\qed




\section{Classification of cubic VNC graphs of order $12p$}
This section is devoted to classifying all connected cubic $\rm VNC$ graphs of order $12p$ for each prime $p$.
The following is the main result of this paper.

\begin{theorem} {\rm \label{thm3}}
A connected cubic graph of order $12p$ for a prime $p$ is a $\rm VNC$ graph if and only if it is isomorphic to
one of $\F084$, $\F204$, ${\rm NC}_{12p}^{i}$ $(0\leq i\leq 9)$.
\end{theorem}

\demo By \cite{CD}, $\F084$ and $\F204$ are connected cubic symmetric graphs of order $12\cdot7$ and $12\cdot 17$,
respectively. By MAGMA \cite{BCP}, $\Aut(\F084)$ and $\Aut(\F204)$ have no regular subgroups.
It follows from Proposition \ref{cayley graph} that $\F084$ and $\F204$ are non-Cayley graphs. By Theorems \ref{thm4}
and \ref{thm2}, the graphs ${\rm NC}_{12p}^{i}$ $(0\leq i\leq 9)$ are connected
cubic $\rm VNC$ graphs of order $12p$ with $p$ a prime.

For the necessity, let $X$ be a connected cubic $\rm VNC$ graph of
order $12p$. By McKay \cite{McKay-Royal}, the numbers of cubic VNC graphs of order $24$, $36$, $60$ and
$84$ are $0,0,8$ and $3$, respectively. It follows that if $p\leq 7$, then we have that $X$
is isomorphic to either ${\rm NC}_{12\cdot 5}^i (0\leq i\leq 6)$, ${\rm NC}_{12\cdot 5}^9$,
${\rm NC}_{12\cdot 7}^i (7\leq i\leq 8)$ or $\F084$.
If $X$ is symmetric and $p>7$, then by Theorem \ref{thm1}, $X\cong \F204$.
In what follows, assume that $p>7$ and that $X$ is non-symmetric.



Let $A=\Aut(X)$. As $X$ is non-symmetric, the vertex-stabilizer $A_v$ is a $2$-group, and
since $X$ is non-Cayley, by Proposition~\ref{cayley graph}, one has
$|A|=2^{m+2}\cdot3\cdot p$ for some positive integer $m$ and $A$ does not contain a regular subgroup.
Suppose that $A$ is non-solvable. Then $A$ has a non-abelian simple composition factor, say $M/N$.
Since $3^{2}\nmid |A|$, by \cite[pp.12--14]{G}, one has $M/N=A_{5}$ or $\PSL(2,7)$ with
$p=5$ or $7$, respectively. This is contrary to our assumption.
Thus, $A$ is solvable. Also, $X\ncong X(3p,2)$ because by the argument preceding the statement of Theorem~\ref{lem1}
we have that $X(3p,2)$ is a Cayley graph. Now it follows from Theorem~\ref{lem1} that $A$
has a normal Sylow $p$-subgroup, say $P$. As $|A|=2^{m+2}\cdot3\cdot p$, one has $P\cong \mz_{p}$.
Let $X_{P}$ be the quotient graph of $X$ relative to $P$,
and let $K$ be the kernel of $A$ acting on $V(X_{P})$. Then $X_{P}$ has order $12$ and
$A/K$ is a vertex-transitive group of automorphisms of $X_P$. Since $P\unlhd A$, the valency of
$X_P$ is $2$ or $3$.

Suppose that $X_P$ is of valency $3$. Then for any vertex $v\in V(X)$,
its neighbors are in three different orbits of $P$. It follows that the stabilizer $K_{v}$ fixes the
neighborhood of $v$ in $X$ pointwise because $K$ fixes each orbit of $P$ setwise.
By the connectedness of $X$, $K_{v}$ fixes each vertex in $V(X)$, forcing $K_{v}=1$. Hence,
$K=PK_{v}=P$ and $A/P\leq\Aut(X_P)$. Noting that $X_P$ is a cubic vertex-transitive graph of order $12$,
by \cite{McKay}, we see that $X_P$ is a Cayley graph and either $|\Aut(X_{P})|=24$, or
$X_{P}\cong X(3,2)$. If $A/P$ contains a subgroup, say $G/P$, acting regularly on $V(X_P)$,
then $G$ acts regularly on $V(X)$, a contradiction. Since $X_P$ is always a Cayley graph, we must have
that $X_{P}\cong X(3,2)$ and $12<|A/P|<|\Aut(X_P)|$. It is easy to see that $\Aut(X(3,2))\cong \mz_{2}^{3}\rtimes S_{3}$.
Consequently, we obtain that $|A/P|=24$, and so $A/P\unlhd \Aut(X_P)$.
Clearly, $O_3(\Aut(X_P))=1$. This implies that
$O_3(A/P)=1$ because $A/P\unlhd\Aut(X_P)$.
In particular, $A/P$ is non-abelian.

Let $C=C_A(P)$. Clearly, $P\leq C$. If $P=C$, then by Proposition \ref{NC}, we have $A/P=A/C\leq \mz_{p-1}$,
forcing that $A/P$ is cyclic, a contradiction. Thus, $P<C$. Let $Q=O_2(C)$. Since $O_3(A/P)=1$, one has $O_3(C/P)=1$.
If $|Q|=2$, then the fact that $O_3(C/P)=1$ implies that $C=P\times Q\cong\mz_{2p}$, and so
$C\leq C_A(C)\leq C_A(P)=C$. It follows that $C=C_A(C)$ and hence
$A/C\leq\Aut(C)\cong\mz_{p-1}$. Since $C/P\cong\mz_2$, it is in the center of $A/P$, and
since $(A/P)/(C/P)(\cong A/C)$ is cyclic, one has $A/P$ is abelian, a contradiction.
So, we must have $|Q|>2$. Clearly, $Q\unlhd C$, and since $C\unlhd A$, one has $Q\unlhd A$.
Recall that $O_2(\Aut(X_P))\cong\mz_2^3$ has three orbits (on $V(X_P)$) each of which has size $2$. Since $A/P\unlhd\Aut(X_P)$,
one has $QP/P\leq O_2(A/P)\leq O_2(\Aut(X_P))$. It follows that
each orbit of $QP/P$ on $V(X_P)$ has size $2$, and so each orbit of $QP$ on $V(X)$ has size $2p$.
Since $Q\unlhd QP$, each orbit of $Q$ on $V(X)$ has size $2$. As $|Q|>2$, the quotient graph $X_Q$ of $X$
relative to $Q$ is a $6p$-cycle, say $X_Q=(B_0, B_1, B_2, \cdots, B_{6p-1})$, and each orbit of $Q$ does not contain an edge.
So, we may assume that $X[B_i,B_{i+1}]\cong 2K_2$ if $i$ is odd and
$X[B_i,B_{i+1}]\cong C_4$ if $i$ is even. Now one may easily see that $X$ is isomorphic to
$X(3p,2)$. This is impossible.

By now, we have shown that $X_P$ is of valency $2$, namely, $X_P\cong C_{12}$.
Without loss of generality,
assume that $V(X_{P})=\{\Delta_{0},\Delta_{1},\cdots,\Delta_{11}\}$, and $\Delta_{i}$
is adjacent to $\Delta_{i+1}$ for each $i\in \mz_{12}$. If $X[\Delta_{i}]$ contains some edges of $X$,
then the connectedness of $X_{P}$ yields $d(X[\Delta_{i}])=1$. This forces that $|\Delta_{i}|=p$ is
even, a contradiction. Thus, $X[\Delta_{i}]$ is an empty graph for every $i\in \mz_{12}$.
Since $X$ is cubic and $d(X_{P})=2$, for any two adjacent orbits $\Delta$, $\Delta^{\prime}$ of $P$,
we have $X[\Delta\cup \Delta^{\prime}]\cong C_{2p}$ or $pK_{2}$. Without loss of generality,
assume that $X[\Delta_{0}\cup \Delta_{11}]\cong pK_{2}$ and $X[\Delta_{0}\cup \Delta_{1}]\cong C_{2p}$.
Then $A/K$ is not edge-transitive on $X_{P}$. By the vertex-transitivity of $A/K$ on $X_{P}$,
we have $A/K\cong D_{12}$. Since $p>3$, the subgroup $K^{\ast}$ of $K$ fixing $\Delta_{0}$ pointwise
also fixes $\Delta_{1}$ and $\Delta_{11}$ pointwise. The connectedness of $X$ gives $K^{\ast}=1$,
and consequently, $K\leq\Aut(\Delta_{0}\cup \Delta_{1})\cong D_{4p}$. Since $K$ fixes
$\Delta_{0}$ setwise, one has $K\cong \mz_{p}$ or $D_{2p}$. Since $|A|>12p$, it follows that
$K\cong D_{2p}$ and hence $|A|=24p$. Since $A/K$ is regular on $V(X_{P})$,
one has $A_{v}=K_{v}\cong \mz_{2}$ and $K=P\rtimes A_{v}$.

Set $C=C_{A}(P)$. Then $P\leq C$. Since $P\unlhd A$, by Proposition \ref{NC},
$A/C\leq \Aut(P)\cong \mz_{p-1}$. As $K=P\rtimes A_{v}\cong D_{2p}$, one has $C_{v}=1$, namely, $C$ is semiregular
on $V(X)$. Then $|C|\ |\ 12p$. Observing that $A/K\cong D_{12}$, it follows that the Sylow $2$-subgroups of
$A$ are non-cyclic. Hence $|C|\neq p,3p$.
If $|C|=2p$ or $4p$, then we could take a Hall $\{2, p\}$-subgroup, say $H$, of $A$ such that $CK\leq H$.
Then $|H|=8p$, and $H/C\unlhd A/C$ since $A/C$ is cyclic. So, $H\unlhd A$, and hence $H/K\unlhd A/K$.
This is impossible because $A/K (\cong D_{12})$ has no normal subgroups of order $4$.
Therefore, $|C|=6p$. Clearly, $C$ has two orbits, say $B_{0}$ and $B_{1}$.
The action of $C$ on each of these two orbits is regular. Furthermore, $d(X[B_{0}])=d(X[B_{1}])=0$, $1$, or $2$.

Clearly, every orbit of $P$ is contained in $B_{0}$ or $B_{1}$. Suppose $X[B_{0}]$ and $X[B_{1}]$ are of valency 0 or 1. Since the induced subgraph of two adjacent orbits of
$P$ is isomorphic to $pK_{2}$ or $C_{2p}$, there are two orbits, say
$\Delta_0$ and $\Delta_1$, of $P$ such that $\Delta_0\subseteq B_0$ and
$\Delta_1\subseteq B_1$ and the subgraph induced by $\Delta_0\cup \Delta_1$
is a $2p$-cycle. Let $M$ be the subgroup of $A$ fixing $\Delta_0\cup \Delta_1$
setwise. Then $M$ acts transitively and faithfully on $\Delta_0\cup \Delta_1$,
so $M\cong D_{4p}$ or $|M|=2p$. Consequently, there is an involution, say $\a$,
in $M$ interchanging $\Delta_0$ and $\Delta_1$. This implies that $C\rtimes\lg \a\rg$
is regular on $V(X)$, a contradiction.

The above argument shows that $d(X[B_{i}])=2$. Recall that for any two adjacent orbits $\Delta$ and $\Delta^{\prime}$ of $P$,
$X[\Delta]\cong X[\Delta^{\prime}]\cong pK_{1}$ and $X[\Delta\cup \Delta^{\prime}]\cong pK_{2}$ or $C_{2p}$. It follows that $X[B_{0}]\cong X[B_{1}]\cong3C_{2p}$. So, $X$ is a $2$-type Bi-Cayley graph over $C$. By viewing $B_0$ and $B_1$ as two copies of $C$, we can let $X={\rm BiCay}(C, R, L, S)$, where $S$ consists of the identity element of $C$, and $|R|=|L|=2$, $R=R^{-1}$,
$L=L^{-1}$ and $C=\lg R\cup L\rg$.

If $C$ is abelian, then $C\cong \mz_{6p}$. By \cite[Theorem~1.1]{ZF2014} and the fact that $p>7$, there is no such 2-type Bi-Cayley but not Cayley graph over C.
Therefore, $C$ is non-abelian. Then $C\cong S_3\times\mz_p$ because $C=C_A(P)$.
 Since $X[B_{0}]\cong X[B_{1}]\cong3C_{2p}$, the fact that $C\cong S_3\times\mz_p$ implies that both $R$ and $L$ consist of an element of order $2p$ and its inverse. For convenience, let $C=\lg a, b\ |\ a^2=b^2=(ab)^3=1\rg\times\lg c\rg\cong S_3\times\mz_p$. By Proposition~\ref{basicprop}, we can take $R=\{ac, ac^{-1}\}$ and $L=\{bc^\ld, bc^{-\ld}\}$. Since $|A|=24p$ and $C\unlhd A$, by Proposition~\ref{one-matching}, $A=C\rtimes\lg\d_\a\rg\cong (S_3\times\mz_p)\rtimes\mz_4$, where $\a\in\Aut(C)$ such that $R^\a=L$ and $L^\a=R$, and $\a$ is not an involution. Since $\a$ interchanges $R$ and $L$, $\a$ must interchange $a$ and $b$ and $\a$ has order $4$. Let $(ac)^\a=bc^i$ with $i=\ld$ or $\ld^{-1}$. Then $(c^2)^\a=((ac)^2)^\a=(bc^i)^2=c^{2i}$. This implies that $i$ is an element of $\mz_p^*$ of order $4$. Thus, $X\cong {\rm NC}_{12p}^9$. \hfill\qed

\Acknowledgements{This work was supported by the National
Natural Science Foundation of China (11671030, 11171020, 11231008) and the Fundamental
Research Funds for the Central Universities (2015JBM110).}

\end{document}